\documentclass[aap,MSNbibl,seceqn,dvips]{arximspdf}

\doi{10.1214/10-AAP694}
\volume{21}
\issue{1}
\pubyear{2011}
\firstpage{266}
\lastpage{282}

\makeatletter

\newtheorem{theorem}{Theorem}[section]

\newtheorem{lemma}[theorem]{Lemma}

\newproclaim{remark}[theorem]{Remark}

\newcommand{\sesq}[2]{\langle #1,#2 \rangle}
\newcommand{\sequiv}{\mathcal{S}}
\newcommand{\timeh}{\bar{T}}
\newcommand{\wienerp}[2]{W^{#1}_{#2}}
\newcommand{\wienerq}[2]{\bar{W}^{#1}_{#2}}
\newcommand{\zcpx}[1]{\tilde{p}_{#1}}
\newcommand{\zcpxd}[1]{p_{#1}}
\newcommand{\prtfpxd}[1]{V_{#1}}
\newcommand{\vol}[2]{\sigma^{#1}_{#2}} %
\newcommand{\dualvol}[2]{\sigma'^{{#1}}_{#2}} %
\newcommand{\drift}[1]{m_{#1}}
\newcommand{\ltrans}[1]{\mathcal{L}_{#1}}
\newcommand{\prtfs}{\mathsf{P}}
\newcommand{\mpr}[2]{\gamma^{#1}_{#2}}
\newcommand{\fwrate}[1]{R_{#1}}
\newcommand{\derprod}[2]{\mathsf{D}^{#1}_{#2}}

\makeatother

\begin{document}
\begin{frontmatter}

\title{Generalized integrands and bond portfolios: Pitfalls and
counter examples}
\runtitle{Generalized integrands and bond portfolios}

\begin{aug}
\author[A]{\fnms{Erik} \snm{Taflin}\thanksref{t1}\corref{}\ead[label=e1]{taflin@eisti.fr}}

\thankstext{t1}{Chair in Mathematical Finance EISTI, and Laboratory
AGM (UMR 8088-CNRS), 95000 Cergy, France.}
\runauthor{E. TAFLIN}
\affiliation{EISTI}
\address[A]{EISTI \\
Avenue du Parc \\
95011 Cergy \\
France \\
\printead{e1}} 
\end{aug}

\received{\smonth{9} \syear{2009}}
\revised{\smonth{3} \syear{2010}}

%
\begin{abstract}
We construct Zero-Coupon Bond markets driven by a cylindrical Brownian
motion in which the notion of generalized portfolio has important
flaws: There exist bounded smooth random variables with generalized
hedging portfolios for which the price of their risky part is $+\infty$
at each time. For these generalized portfolios, sequences of the prices
of the risky part of approximating portfolios can be made to converges
to any given extended real number in $[-\infty, \infty].$
\end{abstract}

\begin{keyword}[class=AMS]
\kwd[Primary ]{60H05}
\kwd[; secondary ]{60G44, 91B28, 91B70}.
\end{keyword}

\begin{keyword}
\kwd{Complete markets}
\kwd{bond markets}
\kwd{generalized integrands}
\kwd{generalized portfolios}
\kwd{replication}.
\end{keyword}

\end{frontmatter}

\section{Introduction}
In this article, we consider continuous time bond markets for which
there exists a unique equivalent martingale measure (e.m.m.).
It is well-known that the uniqueness of the e.m.m. does not in general
imply that such a market is complete. We have here adopted the
standard definition of complete market, which we shall call $L^\infty
$\textit{-completeness} and which reads, omitting details:

\begin{quote}
\textit{Every random variable $X$ in $L^\infty$ is replicable by an
admissible $H'$-valued self-financed portfolio process $\theta,$ where
$H'$ is the dual of $H,$ the state space of the price process.}
\end{quote}

To our knowledge, such \textit{noncompleteness results} were first
established in \cite{BjKaRu97} and \cite{BjMaKaRu97} (see
Proposition 4.7 of \cite{BjKaRu97} and Proposition 6.9 of
\cite{BjMaKaRu97}). The considered price model was a jump-diffusion
model with a finite dimensional Brownian motion (B.m.) and an infinite
mark space, and $H$ was the sup normed Banach space of continuous
functions on $[0,\infty[$ with vanishing limit at $\infty.$ It was also
proved that this market is \textit{approximately complete,} that is,
the subspace of replicable random variables is dense in $L^\infty,$ if
and only if the e.m.m. is unique (see Proposition 6.10 and Theorem 6.11
of \cite{BjMaKaRu97}).

Similar results were proved in \cite{ETBondCompleteness} (see
Theorem 4.1, Theorem 4.2 and Remark~4.6 of
\cite{ETBondCompleteness}) for the case of price models introduced
in \cite{IEET2005}, where the price is a $H$-valued process driven
by a standard cylindrical B.m. (cf. \cite{DaPratoZabczyk}) and where
$H$ is a Sobolev space of continuous functions. Also various
topological vector spaces $A$ for which these markets are $A$-complete
(change $L^\infty$ to $A$ in the above definition) were specified in
Theorem 4.3 of \cite{ETBondCompleteness}. Hedging in the case of a
Markovian price model was considered in \cite{CarmonaTehr}.

The notion of admissible portfolio was weakened in
\cite{DonnoPratelli2005} to that of \textit{generalized self-financed
bond portfolio}, which reestablished, for very general price
processes having a unique e.m.m., the $L^\infty$\textit{-completeness}
of the market, but with $H'$ in the above definition replaced by the
larger set $\mathcal{U}$ of bounded and unbounded linear forms in $H$
(see the discussion in Section 4 of \cite{DonnoPratelli2005} and see
also \cite{DonnoPratelli2004}).

The aim of the present paper is to study and establish properties of
generalized self-financed bond portfolios. In particular, we are
interested in the price of the risky part (or equivalently, in the
price of the risk-free part) of generalized bond portfolios, for which
the separation into risk-free and risky part makes sense. To this end,
simple price models driven by a standard cylindrical B.m., of the kind
introduced in \cite{IEET2005} and \cite{IEET2007} and with
constant volatility operator, are considered. It is proved that the
price model can be chosen such that some generalized self-financed bond
portfolio will have properties to be handled
with care and which can even limit the mathematical and practical
usefulness of generalized portfolios. In fact, (see Theorem \ref{main
th1}):
\begin{longlist}
\item[(a)] There exist bounded smooth random variables, hedgeable in
the sense
of \cite{DonnoPratelli2005} by a unique generalized self-financed
bond portfolio $(x,\mu),$ whose risky part $\mu^1$ is unique and is a
positive $C^{\infty}$ density. The price of $\mu^1$ is $+\infty$ at
each time. Equivalently, it requires to hold a loan of infinite amount,
at each time.

\item[(b)] For all ``admissible'' utility functions, there exists a unique
well-defined optimal wealth $\hat X,$ solution of the optimal expected
utility problem. $\hat X$ is hedgeable in the sense of
\cite{DonnoPratelli2005} by a a unique generalized portfolio
$(x,\mu).$ Also here this generalized portfolio requires to hold a loan
of infinite amount, at each time.

\item[(c)] In each one of the cases (a) and (b), approximate portfolios
converging to $(x,\mu)$ can be chosen such that the sequence of the
prices of their risky part converges to any given extended real number
in $[-\infty, \infty].$
\end{longlist}

Theorem \ref{main th1} gives counter examples to some statements in
\cite{DonnoPratelli2005} (see Remark \ref{rm2}). Results analogous to
those of this paper should apply to other infinite dimensional markets,
as in \cite{Pham2003}.

The present article is a motivation for future research on the hedging
problem in bond markets treated as a super-replication problem under
constraints instead of replication by ``standard'' or generalized
portfolios. \\


\section{Mathematical set-up and the market model}

We shall use a simple case of the Hilbert space Zero-Coupon Bond
models of \cite{IEET2005} and \cite{IEET2007}. The
Zero-Coupon Bond price curves belong to a Hilbert space $H,$ of
continuous functions on $[0,\infty[.$ In this paper, we choose
$H=H^{1}([0,\infty[),$ the Sobolev space of order $1$ of real-valued
functions on $[0,\infty[.$
Let $\ltrans{}$ be the contraction semi-group of left translations in
$L^{2}([0,\infty[),$ let $\partial$ be its infinitesimal generator and
for a positive integer $n \geq0$ let $H^{n}([0,\infty[)$ be the
subspace of functions $f$ such that $[0,\infty[\,\ni a \mapsto
\ltrans{a} f \in L^{2}([0,\infty[)$ is $n$-times continuously
differentiable. $H^{n}([0,\infty[)$ is a Hilbert space for the norm
defined by
\begin{equation} \label{Hn}
\Vert f\Vert_{H^{n}}=\Biggl( \int_{0}^{\infty}
\sum_{i=0}^n|\partial^i f(x)|^{2}\,dx\Biggr)^{1/2}
\end{equation}
and $\ltrans{}$ (restricted to $H^{n}([0,\infty[)$) is a contraction
semi-group in $H^{n}([0,\infty[).$ Pointwise multiplication
$H^{n}([0,\infty[) \times H^{n}([0,\infty[) \ni(f,g) \mapsto fg \in
H^{n}([0,\infty[)$ is continuous for $n \geq1.$

A real-valued    bi-linear form
$\sesq{\cdot}{\cdot},$ where
\begin{equation}\label{sesqprod}
\sesq{f}{g}=\int_0^\infty f(x)g(x)\,dx,
\end{equation}
is first defined for (real) tempered distributions $f$ with support
contained in $[0,\infty[ $ and for (real) tempered test functions $g$
on $\mathbb{R}.$ $H^{-n}([0,\infty[)$ is the subset of all such $f,$
for which the mapping $g \mapsto\sesq{f}{g}$ has a continuous
extension to $H^{n}([0,\infty[).$ The dual $(H^{n}([0,\infty[))'$ of
$H^{n}([0,\infty[)$ is identified with $H^{-n}([0,\infty[)$ and we
write $H'=(H^{1}([0,\infty[))'.$

We consider a time interval $\mathbb{T}=[0,\timeh],$ where $\timeh
>0$ is a finite time-horizon. The random source is an infinite
dimensional $\ell^{2}$-cylindrical Brownian motion
$\wienerp{}{}=(\wienerp{1}{}, \ldots,\wienerp{n}{}, \ldots)$ on a
complete filtered probability space $(\Omega
,P,\mathcal{F},(\mathcal{F}_t)_{t \in\mathbb{T}}),$ where
$\mathcal{F}=\mathcal{F}_{\timeh}$ and the filtration is generated by
the independent Brownian motions $\wienerp{n}{},$ $n \geq1.$

The price at time $t \in\mathbb{T}$ of a Zero-Coupon Bond with time
to maturity $x \geq0$ is denoted $\zcpx{t}(x)$ and the corresponding
discounted price $\zcpxd{t}(x).$ By convention, $\zcpx{t}(0)=1.$

In this paper, we shall use a time independent volatility operator
$\vol{}{}\!\in\! L_2(\ell^{2},H^2),$ the space of Hilbert--Schmidt
operators from $\ell^{2}$ to $H^2([0,\infty[).$ For $z \in\ell^{2},$
$\vol{}{} z= \sum_{i \geq1} \vol{i}{} z^i,$ where the functions
$\vol{i}{}$ satisfy $\vol{i}{}(0)=0.$ Moreover, we impose that
$\vol{i}{} \in C^\infty([0,\infty[),$ that $\vol{i}{}$ has compact
support, that the set $\{\vol{1}{}, \ldots, \vol{i}{}, \ldots\}$ is
linearly independent and total in the subspace of functions $f \in
H^1([0,\infty[)$ satisfying $f(0)=0.$ In particular, it follows that
$\vol{}{}$ is injective.

A drift function $\drift{}$ is given such that $\drift{}=\vol{}{}
\mpr{}{},$ for a time independent market price of risk $\mpr{}{} \in
\ell^2.$ In particular, it follows that $\drift{} \in H^2([0,\infty[),$
since $\vol{}{} \in L_2(\ell^{2},H^2),$ and that $\drift{}(0)=0.$

The discounted price $\zcpxd{}$ is a continuous $H$-valued process
satisfying
\begin{equation}\label{bond dyn integ 1 p*}
\zcpxd{t}= \ltrans{t }\zcpxd{0}+\int_{0}^{t}\ltrans{t-s}(\zcpxd
{s}\drift{})\,ds
+\int_{0}^{t}\ltrans{t-s} (\zcpxd{s} \vol{}{})\,d\wienerp{}{s},
\end{equation}
where $\zcpxd{0} \in H^2([0,\infty[)$ is a strictly positive function
with $\zcpxd{0}(0)=1.$ Here, the notations of pointwise multiplication
are used, so explicitly for the integrand in the second integral:
$(\ltrans{a}(\zcpxd{s} \vol{}{})z)(x)=\sum_{i\geq1}\zcpxd{s}(x+a)
\vol{i}{}(x+a)z^i$ for all $z \in\ell^{2}$ and $a,x \geq0.$

Equation (\ref{bond dyn integ 1 p*}) has a unique $H$-valued mild
solution $\zcpxd{}$ (see \cite{IEET2005} and \cite{IEET2007}
for properties of the solution of (\ref{bond dyn integ 1 p*})). This
solution is a strong solution and it satisfies the following equation
in $H,$ which shows that $\zcpxd{}$ is a $H$-valued semi-martingale:
\begin{equation}\label{bond dyn integ 1 dp*}
d\zcpxd{t}= (\partial\zcpxd{t}+ \zcpxd{t} \drift{})\,dt
+\zcpxd{t} \vol{}{}\,d\wienerp{}{t}.
\end{equation}
For later reference, we note that it follows from Theorem 2.2 of
\cite{IEET2005}, that the mapping $[0,\infty[\,\ni x \mapsto
p(x)$ is a continuous mapping into the space of real semimartingales
$S(P)$ endowed with the semimartingale topology, cf.
\cite{Protter2004}.

A portfolio, also called ``standard portfolio'' in this paper, is an
$H'$-valued progressively measurable process $\theta$ defined on
$\mathbb{T}.$ If $\theta$ is a portfolio, then its discounted value
at time $t$ is
\begin{equation}
\prtfpxd{t}(\theta)=\sesq{\theta_{t}}{ \zcpxd{t}}.
\label{wealth *}
\end{equation}
$\theta$ is an \textit{admissible portfolio} if\setcounter{footnote}{1}\footnote{In this paper,
all considered admissible portfolios will also satisfy
$\prtfpxd{t}(\theta) \geq C$ a.e. $(t,\omega)$ for some $C \in
\mathbb{R}$ depending on $\theta.$}
\begin{eqnarray}\label{prtf norm}
\qquad\Vert \theta\Vert^{2}_{\prtfs}
&=&E\biggl( \int_{0}^{\timeh} (\Vert \theta_{t} \Vert^{2}_{H'}+ \Vert\dualvol
{}{}\theta_{t}\zcpxd{t} \Vert^{2}_{\ell^{2}})\,dt +\biggl(\int_{0}^{\timeh} | \sesq{\theta_{t}}{\zcpxd{t} \drift{}} | \,dt\biggr)^{2}
\biggr) < \infty,
\end{eqnarray}
where $\dualvol{}{}$ is the adjoint of $\vol{}{}$ defined by
$\sesq{f}{\vol{}{}x}=(\dualvol{}{}f,x)_{\ell^{2}},$ for all $f \in H'$
and $x \in\ell^{2}.$ Explicitly, we have:
\begin{equation}\label{adjoint p*}
\dualvol{}{}f=(\sesq{f}{\vol{1}{}},\ldots,\sesq{f}{\vol{i}{}},
\ldots).
\end{equation}
The set of all admissible portfolios defines a Banach space $\prtfs$
for the norm $ \Vert\cdot\Vert_{\prtfs}.$ A portfolio $\theta\in\prtfs$
is by definition \textit{self-financed} if
\begin{equation}
d\prtfpxd{t}(\theta)=\sesq{\theta_{t}}{\zcpxd{t}\drift{}}\,dt
+\sum_{i \in\mathbb{N}^{*}}\sesq{\theta_{t}}{\zcpxd{t}\vol{i}{}}\,d\wienerp{i}{t}.
\label{SFprtf}
\end{equation}

There is a unique e.m.m. (equivalent martingale measure) $Q.$ It is
given by $dQ /dP= \xi_{\timeh},$ where
\[
\xi_t=\exp\bigl((\mpr{}{}, \wienerp{}{t})_{\ell^2} - \tfrac{1}{2}\|
\mpr{}{}\|_{\ell^2}^2 t\bigr).
\]
By Girsanov's theorem the $\wienerq{i}{},$ $i\geq1$, where $\wienerq
{i}{t}=\wienerp{i}{t} + \mpr{i}{} t,$ are independent $Q$-B.m. Obviously,
\begin{equation}
\zcpxd{t}= \ltrans{t }\zcpxd{0}
+\int_{0}^{t}\ltrans{t-s} (\zcpxd{s} \vol{}{})\,d\wienerq{}{s}.
\label{bond dyn integ 1 p}
\end{equation}

We shall only consider derivative products with \textit{discounted
pay-off} belonging to the (Fr\'echet) space $\derprod{}{0},$ which by
definition is the intersection of all the spaces
$L^{p}(\Omega,Q,\mathcal{F}),$ $1 \leq p < \infty.$
Such a derivative $X$ has a unique decomposition as a stochastic
integral w.r.t. $\wienerq{}{}$ (cf. \cite{DaPratoZabczyk} and Lemma
3.2 of \cite{ETBondCompleteness})
\begin{equation} \label{WQ decomp}
X=E_Q[X] +\int_0^{\timeh} (x_t,d\wienerq{}{t}),
\end{equation}
where $x$ is a progressively measurable $\ell^2$-valued process
satisfying
\begin{equation} \label{x space}
x \in L^{p}(\Omega,Q, L^{2}( \mathbb{T},\ell^{2})),\qquad1 \leq p <
\infty.
\end{equation}
It is important to have information about the decay properties of
$x_t^n$ for large $n,$ to study hedging properties of $X.$ We therefore
also introduce the spaces of derivative products $\derprod{}{s},$ $s >
0.$
$\derprod{}{s}$ is the subspace of all $X \in\derprod{}{0}$ such that
the integrand $x$ in (\ref{WQ decomp}) satisfies
\begin{equation} \label{Ds}
\biggl(\int_0^{\timeh} \|x_t\|_{\ell^{s,2}}^2\, dt \biggr)^{1/2} \in\derprod{}{0},
\end{equation}
where $\ell^{s,2} \equiv\ell^{s,2}(q)$ is the Hilbert space of real
sequences endowed with the norm
\begin{equation}\label{ls,2}
\|y\|_{\ell^{s,2}}=\biggl(\sum_{i \in
\mathbb{N}^{*}}q_i^{2s}|y^{i}|^{2}\biggr)^{1/2},
\end{equation}
where $q_i \geq1$ is a given increasing unbounded sequence of real
numbers. (See Remark 4.7 of \cite{ETBondCompleteness}, where
$q_i=(1+i^{2})^{1/2}$ was used.)

Later we shall also impose $X$ to be smooth in the sense of Malliavin.

A \textit{hedging portfolio} $\theta$ of $X$ is a self-financed
portfolio $\theta\in\prtfs$ such that
$\sesq{\theta_{\timeh}}{\zcpxd{\timeh}}=X,$ which then is called
replicable.

Bounded and smooth $X$ are not always replicable, see Remark 4.6 and
Theorem~4.1 of \cite{ETBondCompleteness}. By the definition of
self-financed portfolio, it follows that a portfolio $\theta\in
\prtfs$ is a hedging portfolio of $X$ satisfying (\ref{x space}) iff
$\forall t\in\mathbb{T}$ and $i \geq1$
\begin{equation} \label{hedge prtf}
\sesq{\theta_t}{\zcpxd{t}}=E_Q[X \vert \mathcal{F}_t],\qquad\sesq{\theta
_t}{\zcpxd{t}\vol{i}{}}=x^i_t,
\end{equation}
where $x$ is given by formula (\ref{WQ decomp}). When it exists, the
solution $\theta\in\prtfs$ is unique. In fact, if $\theta$ and
$\phi$
are two solutions, then the second formula in (\ref{hedge prtf}) gives
that $\sesq{\theta_t-\phi_t}{\zcpxd{t}\vol{i}{}}=0,$ for all $i
\geq
1.$ Since the set of the $\vol{i}{}$ is total in the subspace of
functions vanishing at $0$ in $H^{1}([0,\infty[),$ it follows that
$\theta_t-\phi_t=b_t \delta_0$ for some real process $b,$ where
$\sesq{\delta_a}{ f}=f(a),$ for $a \geq0.$ The first formula of
(\ref{hedge prtf}) then gives that
$0=\sesq{\theta_t-\phi_t}{\zcpxd{t}}=b_t \zcpxd{t}(0).$ So $b_t=0,$
since $\zcpxd{t}(0) >0.$

A self-financed discounted risk-free investment, with discounted
value\break
$\prtfpxd{t}(\theta)=1,$ is realized by the portfolio
\begin{equation} \label{bank account}
\theta_t=\exp\biggl(\int_0^t \fwrate{s}(0)\,ds\biggr)\delta_0,
\end{equation}
where the instantaneous forward rate $\fwrate{t}(x)$ at $t$ for time
to maturity $x$ is defined by
\[
\fwrate{t}(x)=-(\partial\ln\zcpxd{t})(x).
\]

In certain cases a portfolio $\theta$ can be separated into a risk-free
part and a risky part. This is the case when $0$ is not in the singular
support of $\theta$ or is an isolated point in the singular support
a.e. $(t,\omega).$ Then $\theta$ has a unique decomposition into a
risk-free part $\psi^0$ and a risky part $\psi^1,$ such that
\begin{equation} \label{decomp bond prtf}
\theta= \psi^0 +\psi^1,\qquad\psi^0_t=a_t \delta_0,\qquad
\operatorname{sing\,supp} \psi^1_t \subset \,]0,\infty[,
\end{equation}
where $a$ is progressively measurable real-valued process. So here
$\sesq{\psi^0_t}{\zcpxd{t}}=a_t \zcpxd{t}(0)$ and
$\sesq{\psi^1_t}{\zcpxd{t}}$ are respectively the discounted risk-free
and risky investments at $t$ corresponding to $\theta.$

The notion of \textit{generalized bond portfolio} was introduced in
ref. \cite{DonnoPratelli2005},
in an attempt to circumvent the problem of the existence of bond
markets with a unique e.m.m., %
but which are not complete in the sense that every sufficiently
integrable r.v. is replicable (by a self-financed admissible bond
portfolio).

Let the product-space $\mathbb{R}^{\mathbb{R}_+}$ be given its natural
product-topology and let $\mathcal U$ be the set of all (bounded and
unbounded) linear forms on $\mathbb{R}^{\mathbb{R}_+}.$ Each element $l
\in\mathcal{U}$ is defined by its domain $\mathcal{D}(l)$ and its
values $l(f)$ for $f \in\mathcal{D}(l).$ Adapted to our mathematical
set-up, a \textit{generalized self-financed bond portfolio} (see
Definition 3.1 of \cite{DonnoPratelli2005}) is a pair $(x,\mu),$
where $x$ is a real number (the value of the generalized portfolio at
$t=0$) and where $\mu$ is a generalized integrand in the sense that
$\mu$ is a $\mathcal{U}$-valued weakly predictable process and there
exist simple integrands $\mu^{(n)},$ that is, $\mu^{(n)} =\sum_i
h^{n,i} \delta_{x_{n,i}}$ where the sum is finite and $ h^{n,i}$ are
bounded predictable real processes, such that
\begin{longlist}
\item[(C$_1$)] $\mu^{(n)}$ converges to $\mu$ a.s. in $\mathcal{U}$
(pointwise),

\item[(C$_2$)] the sequence $Y^n,$ where $Y^n_t=\int_0^{s}
\sum_i \sesq{\mu^{(n)}_{s}}{\zcpxd{t}\vol{i}{}}\,d\wienerq{i}{s}$
converges to a limit process $Y \in S(P).$ $Y_t$ is also denoted
\begin{equation} \label{gen integ}
Y_t= \int_0^{t} \sum_i \sesq{\mu_{s}}{\zcpxd{s}\vol
{i}{}}\,d\wienerq{i}{s}.
\end{equation}
\end{longlist}
The limit $Y$ of $Y^n$ is independent of the sequence $(\mu^{(n)})_{n
\geq1}.$ We recall that, more generally (see Theorem 2.4 of
\cite{DonnoPratelli2005}), if $\mu^{(n)}$ is a sequence of generalized
integrands satisfying (C$_2$) then there exists a generalized integrand
$\mu$ such that equality (\ref{gen integ}) is satisfied.

The discounted value process of the generalized portfolio $(x,\mu)$ is
by definition $x +Y.$ For every $x \in\mathbb{R}$ and portfolio $\mu
\in\prtfs,$ $(x,\mu)$ is a generalized self-financed bond portfolio.
A generalized self-financed bond portfolio $(x,\mu)$ is called
\textit{generalized hedging portfolio} of $X$ when
\begin{equation} \label{eq1 main th1}
X=x+\int_0^{\timeh} \sum_i \sesq{\mu_{t}}{\zcpxd{t}\vol
{i}{}}\,d\wienerq{i}{t}.
\end{equation}

\section{Main results}

A natural question is what are the sequences %
of risk-free and risky investments %
permitting to realize a sequence of approximations $(x,\mu^{(n)}),$
satisfying (C$_1$) and (C$_2$), of a generalized self-financed bond
portfolio $(x,\mu).$ What are the limits of these sequences, if they
exist, and are they independent of the choice of the approximating sequence?
More precisely and generally (cf. Theorem 2.4 of
\cite{DonnoPratelli2005}), let $(\mu^{(n)})_{n \geq1}$ be a sequence
of integrands in $\prtfs$ (i.e., portfolios) satisfying (C$_1$) and
(C$_2$), with the corresponding sequence $(Y^n)_{n \geq1}.$
Self-financed portfolios $\theta^{(n)} \in\prtfs$
are then defined by (cf. Proposition 2.5 of \cite{IEET2005})
\begin{equation} \label{decomp gen bond prtf 1}
\theta^{(n)}= b^n \delta_0 +\mu^{(n)},\qquad b^n_t=
\bigl(x+Y^n_t-\bigl\langle\mu^{(n)}_t,\zcpxd{t}\bigr\rangle\bigr)/ \zcpxd{t}(0).
\end{equation}
If the decomposition (\ref{decomp bond prtf}) applies to the portfolios
$\mu^{(n)},$ with risky part $\mu^{(n)1},$ then it follows that the
self-financed portfolio $\theta^{(n)}$ has a unique decomposition
\begin{eqnarray} \label{decomp gen bond prtf 2}
\theta^{(n)}& =& \theta^{(n)0} +\theta^{(n)1},\nonumber\\[-8pt]\\[-8pt]
\theta^{(n)0}_t&=&a^n_t \delta_0, \theta^{(n)1}=\mu^{(n)1},\qquad
\operatorname{sing\,supp} \theta^{(n)1}_t \subset \,]0,\infty[ .\nonumber
\end{eqnarray}
The real-valued process $a^n,$ which is the investment in the risk-free
asset, is then given by
\begin{equation} \label{decomp gen bond prtf 3}
a^n_t= \bigl(x+Y^n_t-\bigl\langle\mu^{(n)1}_t,\zcpxd{t}\bigr\rangle\bigr)/ \zcpxd{t}(0).
\end{equation}
We will come back later to the above questions concerning the possible
limits of the sequence $a^n_t$ of risk-free investments, by studying
the sequence $r^n_t=\sesq{\mu^{(n)1}_t}{\zcpxd{t}}$ of discounted risky
investments.

Another natural question is what risk-free and risky investments are
required to realize the generalized self-financed portfolio $(x,\mu).$
Suppose that the risky part, lets call it $\mu^1,$ is well-defined.
Then if $\mu^1_t$ is (a.s.) a positive density, its discounted value
$\sesq{\mu^1_t}{\zcpxd{t}} \in[0,\infty]$ (a.s.) is well-defined.
This can easily be generalized to the case where the limit of $\int_0^x
\mu^1_t(y) \zcpxd{t}(y)$ as $x \rightarrow\infty$ makes sense. In
these cases, the risk-free investment is obtained as in (\ref{decomp
gen bond prtf 3})
\[
a_t= (x+Y_t-\sesq{\mu^{1}_t}{\zcpxd{t}})/ \zcpxd{t}(0).
\]

We shall construct a bond market and generalized self-financed bond
portfolios $(x,\mu),$ whose realization require an infinite short
position in the risk-free asset (i.e., loan) at each instant $t \in
\mathbb{T}.$ More precisely, to have a clear separation between the
investment into the risk-free asset and the risky assets, we construct
a market
and generalized portfolios $(x,\mu)$ satisfying:
\begin{longlist}
\item[(P$_1$)]
$\nu$ is an element in $\mathcal{U}$ with domain ($\operatorname{ls}$
denotes linear span)
\begin{equation} \label{P0}
\mathcal{D}(\nu) =\operatorname{ls}\bigl(C^{\infty}_0(]0, \infty[) \cup \{
\zcpxd{0} \}\bigr).
\end{equation}
The restriction of $\nu$ to $C^{\infty}_0(]0, \infty[)$ is a function
$\nu^1 \in C^{\infty}([0,\infty[),$ $\operatorname{supp} \nu^1 \subset
[3/4,\infty[ ,$ $\sesq{\nu}{\zcpxd{0}}=0$ and
\begin{equation} \label{P01}
\lim_{x \rightarrow\infty} \int_0^x \nu^1(y) \zcpxd{0}(y)\,dy =
\infty.
\end{equation}
$\mu_t \in\mathcal{U}$ a.s. has domain
\begin{equation} \label{P1.1}
\mathcal{D}(\mu_t)=\operatorname{ls}\bigl(C^{\infty}_0(]0, \infty[) \cup
\{\zcpxd{t}\} \bigr)
\end{equation}
and
\begin{equation} \label{P1}
\sesq{\mu_t}{f}=\alpha_t \biggl\langle\nu,f\frac{\zcpxd{0}}{\zcpxd
{t}}\biggr\rangle,\qquad f
\in \mathcal{D}(\mu_t),
\end{equation}
where $\alpha$ is a strictly positive continuous adapted uniformly
bounded (in $t$ and $\omega$) process.
The discounted total risky investment is
\begin{equation} \label{P1.2}
\lim_{x \rightarrow\infty} \int_0^x \mu^1_t(y) \zcpxd{t}(y)\,dy =
\infty \qquad \mbox{a.e. }(t,\omega).
\end{equation}
\item[(P$_2$)] For all $C \in[-\infty, \infty],$ $\mu$ is the limit
in the sense of (C$_1$) and (C$_2$) of a sequence
$(\mu^{(n)})_{n \geq1}$ of continuous linear functionals on $H$ such
that, a.s.
\begin{eqnarray} \label{P2.0}
\mu^{(n)}&=&\mu^{(n)0}+\mu^{(n)1},\qquad \forall t \in\mathbb{T} \quad\operatorname{supp} \mu^{(n)0}_t \subset
\{0\},\nonumber\\[-8pt]\\[-8pt]
 \mu^{(n)1}_t &\in& C^{\infty
}_0(]1/2, \infty[)\nonumber
\end{eqnarray}
and
\begin{equation} \label{P2}
\lim_{n \rightarrow\infty} a^n_0 = C ,
\end{equation}
where $a^n$ is defined by (\ref{decomp gen bond prtf 3}). Moreover, if $C= -\infty$ (resp. $C$ is finite and $C=+\infty$) then
\begin{equation} \label{P2.1}
\forall t \in\mathbb{T}, \qquad\lim_{n \rightarrow\infty} a^n_t =
-\infty \qquad\mbox{(resp. finite and }+\infty).
\end{equation}
\end{longlist}

\begin{remark} \label{rm1}
The definition of $\nu$ makes sense since $\zcpxd{0}$ is not in $K=\{f
\in H \mid f(0)=0\},$ the closure of $C^{\infty}_0(]0, \infty[)$ in
$H.$ The formula (\ref{P1}) makes sense since, according to Theorem
21 of \cite{IEET2007}, $\| \ell_t / \zcpxd{t} \|_{L^{\infty}} <
\infty,$ where $\ell_t = \mathcal{L}_t \zcpxd{0}.$ So $f \ell_t /
\zcpxd{t} \in\mathcal{D}_0$ a.s.
\end{remark}

An admissible utility function $U$ is (in this article), a strictly
concave and increasing $C^2$ function on $]0,\infty[$ satisfying
conditions, stated in (\ref{condition I}), strengthening the Inada
conditions.
Let $\mathrm{I}$ be the inverse of $U'$ and assume that there exists $C,p
>0$ such that
\begin{eqnarray} \label{condition I}
U'(]0,\infty[ )\,&=&\, ]0,\infty[\quad\mbox{and}\nonumber\\[-8pt]\\[-8pt]
 | \mathrm{I}(x) |+| x \mathrm{I}'(x) | &\leq& C (x^{p}+x^{-p}),\qquad
 x>0.\nonumber
\end{eqnarray}

We shall consider the optimal portfolio problem. For an admissible
utility function $U$ and an initial investment of $E_Q[\mathrm{I}(y
\xi_{\timeh})],$ for given $y>0,$ the optimal final wealth is given by
\begin{equation} \label{opt X}
\hat{X}=\mathrm{I}(y \xi_{\timeh})
\end{equation}
and $\hat{X} \in L^q,$ for all $1 \leq q < \infty,$ cf. Theorem 3.3 of
\cite{IEET2005}.

We can now state the main results (in which risky means that $0$ is not
in the singular support).

\begin{theorem} \label{main th1}
One can choose an initial condition $\zcpxd{0},$ a time-independent
volatility operator $\vol{}{}$ and a time-independent drift function
$\drift{}$ such that:
\begin{longlist}
\item[A.] The $\vol{i}{} \in C^{\infty}_0(]0, \infty[),$ $\vol{}{}
\in
L_2(\ell^{2},H^2([0,\infty[))$ is injective and $\zcpxd{0}(x)=e^{-a
x},$ for some $a>0.$ The drift $\drift{} \in H^2([0,\infty[)$ and the
market price of risk $\mpr{}{} \in\ell^2.$

\item[B.] For all admissible utility functions $U$ and $y >0,$ $\hat
{X}$ given by (\textup{\ref{opt X}}) has a generalized hedging portfolio
$(E_Q[\hat{X}],\mu)$ satisfying (\textup{\ref{eq1 main th1}}) and with the
properties %
\textup{(P}$_1)$ and \textup{(P}$_2)$. The risky part of $(E_Q[\hat{X}],\mu)$ is unique.

\item[C.] There exists a bounded smooth r.v. $X$ having a generalized
hedging portfolio $(E_Q[X],\mu)$ satisfying (\textup{\ref{eq1 main th1}}) and
satisfying \textup{(P}$_1)$ and \textup{(P}$_2)$ with $\nu^1$ positive. The risky part of
$(E_Q[X],\mu)$ is unique and positive.
\end{longlist}
\end{theorem}

\begin{remark} \label{rm2} If $(x,\mu)$ is a generalized hedging
portfolio given by Theorem~\ref{main th1} then:
\begin{longlist}
\item[1.] Since (P$_1$) is satisfied it follows that the value of the risky
part of $(x,\mu),$ is infinite and that the
realization of $(x,\mu)$ requires an infinite short position in the
risk-free asset (i.e., loan) at each instant $t \in\mathbb{T}.$
\item[2.] According to (P$_2$), the sequence of prices, at $t=0,$ $r^n_0$ of
the risky part (or $a^n_0$ of the risk-free part) of approximating
sequences $(x,\mu^{(n)})$ give no information concerning the value of
the risky part (or of the risk-free part) of $(x,\mu).$
As matter of fact for the given $(x,\mu),$ one can choose an
approximating sequence $(x,\mu^{(n)})$ such that the limit of $a^n_0$
is equal to any extended real number in $[-\infty, \infty].$
\item[3.] $\zcpxd{t} \in\mathcal{D}(\mu_t)$ [in fact
$\sesq{\mu_t}{\zcpxd{t}}=0$ a.e., according to (P$_1$)], which is a
condition in a discussion in \cite{DonnoPratelli2005} (second
paragraph after Definition 3.1). The preceding points 1 and 2 of this
remark are counter examples the conclusions of that discussion.
\end{longlist}
\end{remark}

\section{Proofs}
Following \cite{ETBondCompleteness}, we introduce for $t \in
\mathbb{T},$ the operator $B_t=\ell_t \vol{}{} \in L_2(\ell^{2},H),$
where $\ell_t = \mathcal{L}_t \zcpxd{0}.$ Here, $B_t$ is deterministic.
Let $B_{t}^{*}$ be the adjoint of $B_{t}$ with respect to the scalar
product $(\cdot,\cdot)_H$ in $H.$ We also introduce
\begin{equation}\label{At}
A_{t}=B_{t}^{*}B_{t},
\end{equation}
which is a strictly positive self-adjoint trace-class operator in
$\ell^{2}.$ We shall impose the following condition [to be verified
after (\ref{ki def})] on the operators $A_t$: There exists $s >0$ and
$k >0$ such that for all $t \in\mathbb{T}$ and $x \in\ell^{2},$
\begin{equation}\label{uniform cond sigma eq}
\|x\|_{\ell^{2}}
\leq k \|(A_{t})^{1/2}x \|_{\ell^{s,2}}.
\end{equation}
When this condition is satisfied, the contingent claims in
$\derprod{}{s}$ are replicable by self-financed portfolios in $\prtfs$
(Theorem 4.3 of \cite{ETBondCompleteness}).

Let $\sequiv$ be the canonical isomorphism of $H$ onto $H'$ defined by
\begin{equation}
\forall f,g \in H,\qquad (f,g)_{H}=\sesq{\sequiv f}{g}
\label{S equiv}
\end{equation}
and let $S_{t}$ be the isometric embedding of $\ell^{2}$ into $H$ equal
to the closure of $B_{t}(A_{t})^{-1/2}.$ We note that if $f \in H$ is
$C^2,$ then
\begin{equation} \label{S equiv 1}
\sequiv f=f-\partial^2 f - (\partial f)(0) \delta_0.
\end{equation}

If $X \in\derprod{s}{},$ with $s>0$ as in (\ref{uniform cond sigma
eq}), then the equations (\ref{hedge prtf}) have a unique solution
$\theta\in\prtfs$ and $\theta= \theta^0 + \theta^1,$ $\theta^0,
\theta^1 \in\prtfs,$ where
\begin{equation}
\theta^{1}_{t}=(l_{t}/\zcpxd{t}) \sequiv\eta_{t},\qquad \eta_{t}(\omega
)=S_{t}(A_{t})^{-1/2}x_{t}(\omega)
\label{prtf theta1}
\end{equation}
and
\begin{equation}
\theta^{0}_{t}=b_{t} \delta_{0},\qquad b_t = (E_Q[X \vert \mathcal{F}_t]-\sesq
{\theta^{1}_{t}}{\zcpxd{t}})/ \zcpxd{t}(0).
\label{prtf theta0}
\end{equation}
For such $X$,
\begin{equation} \label{thta1 pt}
\sesq{\theta^{1}_{t}}{\zcpxd{t}}=(S_{t}(A_{t})^{-1/2}x_{t},l_{t})_H .
\end{equation}

We shall now construct the volatility operator $\vol{}{}$ and drift
function $\drift{}$ of Theorem \ref{main th1}. For a given $a>0$, we
define $\zcpxd{0}$ by
\begin{equation} \label{def p0}
\zcpxd{0}(x)=\exp(-a x).
\end{equation}
$C^{\infty}_0(]0, \infty[)$ is dense in the (closed) subspace $K$ of
functions $f \in H,$ satisfying $f(0)=0.$ Let $h_1 \in C^{\infty}_0(]0,
\infty[)$ be such that $h_1 \geq0,$ $\operatorname{supp} h_1 \subset
[3/4,5/4]$ and $\|h_1\|_H=1.$ $h_{2n-1} \in C^{\infty}_0(]0, \infty[),$
$n>1$ is defined by $h_{2n-1}(x)=h_1(x-2n+2)$ if $x \geq2n-2$ and
$h_n(x)=0$ if $0 \leq x <2n-2.$ The set of functions $\{h_{2n-1}\}_{n
\geq1}$ is orthonormal in $K$ and
\[
\operatorname{supp} h_{n} \subset \bigl[n - \tfrac{1}{4},n + \tfrac{1}{4}\bigr] ,\qquad n
\mbox{ odd}.
\]
We complete it to an orthonormal basis $\{h_i\}_{i=1}^{\infty} \subset
C^{\infty}_0(]0, \infty[)$ of $K.$
Then \mbox{$h_i/\zcpxd{0} \in K.$}

Let the volatility functions satisfy
\begin{equation} \label{def vol}
\vol{i}{}=k_i h_i / \zcpxd{0},\qquad
k_i \neq0,\qquad \mbox{s.t.}\quad
\sum_{i \geq1} i^2 k_i^2 (1+\|h_i/ \zcpxd{0} \|_{H^2}^2) < \infty.
\end{equation}
The conditions $\vol{}{} \in L_2(\ell^{2},H^2)$ and $\vol{i}{}(0)=0$
are then satisfied and the set $\{ \vol{}i{}\}_{i=1}^{\infty} $ is by
construction linearly independent and total in $K.$

The definition of $B_t$ gives
\begin{equation} \label{B eq1}
B_ty=e^{-at} \sum_{i \geq1} k_i h_i y_i\quad\mbox{and}\quad
(B_{t}^{*}f)^i=e^{-at} k_i (h_i,f)_H.
\end{equation}
It follows that
\begin{equation} \label{At 1}
(A_ty)^i=e^{-2at}k_i^2 y_i \quad\mbox{and}\quad (A_t^{1/2}y)^i=e^{-at} |k_i| y_i.
\end{equation}

It then follows that $(A_t)^{-1/2}$ and $(A_0)^{-1/2}$ have the same
domain and that after closure
\begin{equation} \label{S isom}
S_ty= \sum_i \operatorname{sgn} (k_i) h_i y_i, \qquad y \in\ell^2.
\end{equation}
So for $y \in\mathcal{D}((A_0)^{-1/2})$
\begin{eqnarray} \label{SA}
S_t(A_t)^{-1/2}y&=&e^{at} \sum_{i \geq1} \frac{1}{k_i} h_i y_i\quad
\mbox{and}\nonumber\\[-8pt]\\[-8pt]
\|S_t(A_t)^{-1/2}y \|_H^2&=&e^{2at} \sum_{i \geq1}
\biggl(\frac{y_i}{k_i}\biggr)^2.\nonumber
\end{eqnarray}
This gives for $y \in\mathcal{D}((A_0)^{-1/2})$:
\begin{equation} \label{eq 1}
(S_t(A_t)^{-1/2}y,\ell_t)_H=(S_0(A_0)^{-1/2}y,\zcpxd{0})_H= \sum_{i
\geq1} \frac{1}{k_i} (h_i,\zcpxd{0})_H y_i.
\end{equation}

We define
\begin{equation} \label{m}
m=\vol{}{} \mpr{}{},\qquad \mpr{i}{}=\frac{1}{i}\qquad\mbox{if $i$
odd}\quad\mbox{and}\quad \mpr{i}{}=0\qquad \mbox{if $i$ even}.
\end{equation}

In the following lemma, $H_c$ stands for the complex linear Hilbert
space $H^1([0,\infty[ , \mathbb{C}).$ The function $[0, \infty[\,
\ni x \mapsto e^{-ax}$ is in $H_c$ for $\Re a >0.$

\begin{lemma} \label{lm numb of zeros}
For every $i$ the function $a \mapsto(h_i, e^{-a \cdot})_{H_c},$ $\Re
a >0,$ extends to an entire analytic function on $\mathbb{C}.$
There is only a countable number of $a \in\mathbb{C}$ such that
\begin{equation} \label{eq lm 1}
(h_i, e^{-a \cdot})_{H_c} = 0\qquad\mbox{for some } i \geq1.
\end{equation}
\end{lemma}

\begin{pf}
For $\Re a > 0,$ $F(a) \in H_c,$ where
$(F(a))(x)=e^{-ax}.$ With an obvious extension of $\sesq{\cdot}{\cdot}$ and
recalling that $h$ is real-valued, we have
\[
\lambda_i(a) \equiv(h_i, F(a))_{H_c}=\sesq{\sequiv h_i}{F(a)}.
\]
According to (\ref{S equiv 1}), the distribution $\sequiv h_i$ has
compact support. The Fourier--Laplace transformation $\lambda_i$ of
$\sequiv h_i$ therefore defines an entire analytic function in
$\mathbb{C}.$
Since $\sequiv h_i \neq0$, the set of zeros $A_i$ of the function
$\lambda_i$ in $\mathbb{C}$ is countable. The set
\[
A= \bigcup_{i \geq1} A_i
\]
is then countable, since it is a countable union of countable sets. $A$
is the set of $a$ that satisfies (\ref{eq lm 1}).
\end{pf}

\begin{pf*}{Proof of Theorem \ref{main th1}}
Obviously,
$\zcpxd{0}$ and the $\vol{i}{}$ are as stated in the theorem. By
construction and $\mpr{}{} \in\ell^2$ according to (\ref{m}), so
statement A is true.

To prove the statement B, we follow Remark 4.6 of
\cite{ETBondCompleteness}. Since
\[
\xi_t=\exp\bigl((\mpr{}{}, \wienerq{}{t})_{\ell^2} + \tfrac{1}{2}\|
\mpr{}{}\|_{\ell^2}^2 t\bigr),
\]
$\hat{X}$ has the representation
\[
\hat{X}=E_{Q}[\hat{X}]+\int_{0}^{\timeh}
E_{Q}[y \xi_{\timeh} \mathrm{I}'(y \xi_{\timeh}) \vert \mathcal{F}_{t}]
\sum_{i \geq1} (-\mpr{i}{})\,d\wienerq{i}{t}.
\]
Let $c=1/\|\mpr{}{}\|_{\ell^2}$ [see (\ref{m})] and $e=c \mpr{}{},$
$Z=\sum_{i \geq1}e^{i} \wienerq{i}{\timeh}.$ This proof is based on
the fact that $e \notin\mathcal{D}((A_0)^{-1/2}),$ according to
(\ref{def vol}) and (\ref{At 1}). The real-valued function $g$ is
defined by
\[
g(z)=-\dfrac{1}{c} h(z)I'(h(z)),\qquad h(z)=y \exp\bigl(z/c + \timeh/(2c^2)\bigr),\qquad z
>0,
\]
$g$ is strictly positive on $]0,\infty[.$ Then
\begin{equation} \label{X decomp}
\qquad\hat{X}=E_{Q}[\hat{X}]+\int_{0}^{\timeh}\sum_{i \geq
1}x^{i}_{t}\,d\wienerq{i}{t},\qquad x_{t}= \alpha_t e,\qquad \alpha_t=E_{Q}[g(Z)
\vert
\mathcal{F}_{t}],
\end{equation}
where the r.v. $\alpha_t$ is strictly positive.

The unbounded linear functional $\nu\in\mathcal{U}$ is defined by its
domain given by formula (\ref{P0}) and by
\begin{equation}\label{nu}
\sesq{\nu}{f} =\biggl\langle\sum_{i \geq1 } \frac{e^i}{k_i} \sequiv
h_i,f\biggr\rangle\qquad \mbox{if } f \in C^{\infty}_0(]0, \infty[) \mbox{ and } \sesq
{\nu}{\zcpxd{0}}=0.
\end{equation}
This definition makes sense, since for given $f$ the sum has only a
finite number of terms and since $\zcpxd{0} \notin K,$ the closure of
$C^{\infty}_0(]0, \infty[)$ in $H.$ We define $\nu^1$ by the sum
\begin{equation} \label{nu prime}
\nu^1(x)=\sum_{i \geq1} \frac{e^i}{k_i} \sequiv h_i (x),\qquad x \geq0.
\end{equation}
Here, at most one term is nonvanishing and it must be a term with an
odd index~$i.$ Due to the properties of $h_i$ for odd $i$ and (\ref{S
equiv 1}), we have $\nu^1 \in C^{\infty}([0,\infty[)$ and $\operatorname{supp}
\nu^1 \subset [3/4,\infty[ .$ Obviously, $\nu^1$ is the
restriction of $\nu$ to $C^{\infty}_0(]0, \infty[).$

In order to construct a generalized self-financed bond portfolio
$(E_{Q}[\hat{X}],\mu),$ with value process $Y,$ where
$Y_t=E_{Q}[\hat{X} \vert \mathcal{F}_{t}],$ we define $\mu$ a.e. $dt
\times dP$ by formulas (\ref{P1.1}) and (\ref{P1}) and with $\alpha$
given by (\ref{X decomp}). This makes sense since $f \mapsto
f\frac{\zcpxd{0}}{\zcpxd{t}}$ maps $\mathcal{D}({\mu_t})$ into
$\mathcal{D}({\nu}).$ Property $(P_1)$ is then satisfied.

The sequence $\{e^{(n)}\}_{n \geq1}$ in $\ell^2$ is defined by
$(e^{(n)})^i=e^i$ for $1 \leq i \leq n$ and $(e^{(n)})^i=0$ for $i >n.$ Let
\[
X^n=E_{Q}[\hat{X}]+\int_{0}^{\timeh} \alpha_t \sum_{i \geq1}
\bigl(e^{(n)}\bigr)^i\, d\wienerq{i}{t},\qquad Y^n_t= E_{Q}[X^n \vert \mathcal{F}_{t}].
\]
As $e^{(n)}$ belongs to the domain of $(A_t)^{-1/2}$ we can proceed
as in Remark 4.8 and Theorem 4.3 of \cite{ETBondCompleteness} to
construct the unique hedging portfolio $\theta^{(n)} = \theta^{(n)0}
+\theta^{(n)1},$ where $\theta^{(n)0}, \theta^{(n)1} \in\prtfs$ are
given by (\ref{prtf theta0}) and (\ref{prtf theta1}). Applying
(\ref{SA}) and (\ref{eq 1}), we obtain
\begin{equation} \label{mu n0}
\qquad\theta^{(n)0}_t =a^n_t \delta_0,\qquad a^n_t=\biggl(E_{Q}[X^n \bigl\vert \mathcal
{F}_{t}]- \alpha_t \sum_{1 \leq i \leq n}\frac{e^i}{k_i} (h_i,\zcpxd
{0})_H\biggr) \big/ \zcpxd{t}(0)
\end{equation}
and
\begin{equation} \label{mu n1}
\theta^{(n)1}_t =\frac{\zcpxd{0}}{\zcpxd{t}} \alpha_t \sum_{1
\leq i \leq n} \frac{e^i}{k_i} \sequiv h_i .
\end{equation}
The sequence $\nu^{(n)} \in H'$ is defined by $\sesq{\nu^{(n)}}{f},$ $f
\in H,$ where
\begin{equation} \label{nu n}
\quad \bigl\langle\nu^{(n)},f\bigr\rangle =\biggl\langle\sum_{1 \leq i \leq n} \frac{e^i}{k_i}
\sequiv h_i,f\biggr\rangle, \qquad\mbox{if } f(0)=0 \mbox{ and } \bigl\langle\nu
^{(n)},\zcpxd{0}\bigr\rangle=0.
\end{equation}

We note that $\nu^{(n)}$ converges to $\nu$ in $\mathcal{U}$:
\begin{equation} \label{nu convergence}
\forall f \in \mathcal{D}(\nu),\qquad \lim_{n \rightarrow\infty} \bigl\langle\nu^{(n)},f\bigr\rangle= \sesq{\nu}{f}.
\end{equation}
Let $\nu^{(n)1}$ be the restriction of $\nu^{(n)}$ to $C^{\infty}_0(]0,
\infty[).$
Due to the properties of $h_i$ for odd $i$ and (\ref{S equiv 1}), %
$\nu^{(n)1} \in C^{\infty}([0,\infty[)$ has compact support,
\[
\operatorname{supp} \nu^{(n)}, \operatorname{supp} \nu^{(n+1)} \subset \bigl[3/4,n +
\tfrac{1}{4}\bigr[ ,\qquad n \mbox{ odd}.
\]

We have the decomposition $\nu^{(n)} = \nu^{(n)0} +\nu^{(n)1},$ where
$\nu^{(n)0}, \nu^{(n)1} \in\prtfs$ are given by
\begin{equation} \label{nu n decomp}
\qquad\quad\nu^{(n)1} = \sum_{1 \leq i \leq n} \frac{e^i}{k_i} \sequiv h_i,\qquad \nu
^{(n)0}_t =b^n \delta_0,\qquad b^n=- \sum_{1 \leq i \leq n} \frac
{e^i}{k_i} (h_i,\zcpxd{0})_H.
\end{equation}

We define $\mu^{(n)}_t \in H'$ a.e. $(t,\omega)$ by
\begin{equation} \label{mu n def}
\mu^{(n)}_t=\alpha_t \frac{\zcpxd{0}}{\zcpxd{t}} \nu^{(n)}.
\end{equation}
We have the decomposition $\mu^{(n)} = \mu^{(n)0} +\mu^{(n)1},$ where
\begin{equation}
\mu^{(n)0}_t=\alpha_t \frac{1}{\zcpxd{t}(0)}\nu^{(n)0},\qquad \mu
^{(n)1}_t=\alpha_t \frac{\zcpxd{0}}{\zcpxd{t}}\nu^{(n)1}=\theta^{(n)1}.
\end{equation}

It follows from formulas (\ref{P1.1}), (\ref{P1}) and (\ref{mu n def})
that $\mu^{(n)}_t$ converges a.e. $(t,\omega)$ to $\mu_t$ in
$\mathcal{U},$ so (C$_1$) is satisfied. Since
\[
\bigl\langle\mu^{(n)}_t,\zcpxd{t}\vol{i}{}\bigr\rangle=\alpha_t \bigl(e^{(n)}\bigr)^i,
\]
it follows that $Y^n$ converges to $Y$ in the topology of square
integrable martingales, which is stronger than the semi-martingale
topology. So also (C$_2$) is satisfied. Therefore, $(E_{Q}[X],\mu)$ is
a generalized hedging-portfolio of $X.$

We now fix $a$ and the $k_i.$ $a>0$ is chosen such that $\lambda_i(a)
\equiv(h_i, e^{-a \cdot})_H \neq0$ for all $i \geq1,$ which is
possible according to Lemma \ref{lm numb of zeros}. Let
\begin{eqnarray} \label{ki def}
\operatorname{sgn}(k_i)&=&\operatorname{sgn}(\lambda_i(a)),\nonumber\\[-8pt]\\[-8pt]
0 &<& |k_i|\leq \min\bigl(|\lambda_i(a)|, (1+\|h_i/ \zcpxd{0} \|
_{H^2}^2)^{-1/2}\bigr)/i^2.\nonumber
\end{eqnarray}
The condition in (\ref{def vol}) is then satisfied.

The sequence
$E_{Q}[X^n \vert \mathcal{F}_{t}]$ converges to $E_{Q}[X \vert
\mathcal{F}_{t}]$ in $L^2(\Omega,Q)$ as $n \rightarrow\infty.$ We have
\[
\mathop{\sum_{1 \leq i \leq n}}_{ i\  \mathrm{odd}}\frac{1}{ik_i}
(h_i,\zcpxd{0})_H \geq
\mathop{\sum_{1 \leq i \leq n}}_{ i\  \mathrm{odd}} i .
\]
The last sum goes to $+\infty$ when $n \rightarrow\infty.$ Since $c>0$
and $\alpha_t >0$ a.s. it follows from (\ref{thta1 pt}) that (\ref{P2})
is satisfied in the case of $C=-\infty.$

We shall impose supplementary conditions on the $k_i$ to ensure that
(\ref{P2}) is satisfied also for $C$ finite and $C=+\infty.$ Let
$J\dvtx
\mathbb{N}^* \rightarrow2\mathbb{N}^*+1$ be defined by
\[
J(n)=n+2 \qquad\mbox{if $n$ is odd}\quad\mbox{and}\quad J(n)=n+1 \qquad\mbox{if $n$ is even}.
\]
For $n$ odd let $d^{(n)} \in\mathbb{R},$ for $n$ even let
$d^{(n)}=d^{(n-1)}$ and define for $n \in\mathbb{N}^*$
\[
\tilde{\nu}^{(n)1}= \nu^{(n)1} +d^{(n)} \sequiv h_{J(n)}.
\]
We define $\tilde{\nu}^{(n)} \in H'$ by
\[
\bigl\langle\tilde{\nu}^{(n)},f\bigr\rangle= \bigl\langle\tilde{\nu}^{(n)1},f\bigr\rangle\qquad \mbox{for }
f \in K \mbox{ and } \bigl\langle\tilde{\nu}^{(n)},\zcpxd{0}\bigr\rangle=0.
\]
$\tilde{\nu}^{(n)}$ converges to $\nu$ in $\mathcal{U}$:
\begin{equation} \label{nu tilde convergence}
\forall f \in \mathcal{D}(\nu),\qquad \lim_{n \rightarrow\infty} \bigl\langle\tilde{\nu}^{(n)},f\bigr\rangle= \sesq{\nu}{f}.
\end{equation}
Since $\sesq{d^{(n)} \sequiv h_{J(n)}}{\zcpxd{0} \vol{j}{}}=d^{(n)} k_j
\delta_{j J(n)}$, it follows that
\begin{equation}\label{eq1 cond k d}
\sum_{j=1}^{\infty}\bigl(\bigl\langle d^{(n)} \sequiv h_{J(n)},\zcpxd{0} \vol
{j}{}\bigr\rangle\bigr)^2=\bigl(d^{(n)}\bigr)^2 \bigl(k_{J(n)}\bigr)^2.
\end{equation}
We impose the following condition, which we for the moment suppose is
possible:
\begin{equation}\label{eq2 cond k d}
\lim_{n \rightarrow\infty} d^{(n)} k_{J(n)}=0.
\end{equation}
$\tilde{\mu}^{(n)}_t$ is defined as in (\ref{mu n def}), but with
$\tilde{\nu}$ instead of $\nu.$ Formulas (\ref{nu tilde convergence})
and (\ref{eq2 cond k d}) imply that $(x,\tilde{\mu}^{(n)})$ is an
approximating sequence for the generalized portfolio $(x,\mu).$

We note that
\[
\bigl\langle d^{(n)} \sequiv h_{J(n)},\zcpxd{0}\bigr\rangle=d^{(n)} e^{-aJ(n)}
(h_1,\zcpxd{0})_{H},
\]
which gives
\[
\bigl\langle\tilde{\nu}^{(n)1},\zcpxd{0}\bigr\rangle= \sum_{1 \leq i \leq n}\frac
{e^i}{k_i} (h_i,\zcpxd{0})_H + d^{(n)} e^{-aJ(n)} (h_1,\zcpxd{0})_{H}.
\]
Similarly as in (\ref{mu n0}), we introduce [recalling that
$\zcpxd{0}(0)=1$] $ \tilde{a}^n_0=E_{Q}[\hat{X}]- \alpha_0
\sesq{\tilde{\nu}^{(n)1}}{\zcpxd{0}},$ which gives
\[
\tilde{a}^n_0= E_{Q}[\hat{X}]- \alpha_0 \biggl( \sum_{1 \leq i \leq
n}\frac{e^i}{k_i} (h_i,\zcpxd{0})_H + d^{(n)} e^{-aJ(n)} (h_1,\zcpxd
{0})_{H}\biggr).
\]
For given $\tilde{a}^n_0,$ this is an equation for $d^{(n)}.$ We now
define for $n \geq1$:
\begin{eqnarray*}
\tilde{a}^n_0&=&C \qquad\mbox{if $C$ is finite} \quad\mbox{and}\\
\tilde{a}^n_0&=&
E_{Q}[\hat{X}]+ \alpha_0 \sum_{1 \leq i \leq n}\frac
{e^i}{k_i}(h_i,\zcpxd{0})_H \qquad\mbox{if $C=\infty$}.
\end{eqnarray*}
In both cases, $\lim_{n \rightarrow\infty} \tilde{a}^n_0 =C.$ For $C$
finite, $d^{(n)}$ is then given by
\[
d^{(n)} =\biggl( E_{Q}[\hat{X}]-C - \alpha_0 \sum_{1 \leq i \leq
n}\frac{e^i}{k_i} (h_i,\zcpxd{0})_H \biggr) \frac{e^{aJ(n)}}{\alpha
_0 (h_1,\zcpxd{0})_{H}}
\]
and for $C=+\infty$ by
\[
d^{(n)} =-2\frac{e^{aJ(n)}}{(h_1,\zcpxd{0})_{H}} \sum_{1 \leq i \leq
n}\frac{e^i}{k_i} (h_i,\zcpxd{0})_H .
\]
The property $d^{(n)}=d^{(n+1)}$ is then satisfied for $n$ odd. For odd
$n,$ we choose $|k_{n+2}|>0$ sufficiently small so that $|d^{(n)}
k_{n+2}| \leq1/n.$ Condition (\ref{eq2 cond k d}) is then satisfied.
This proves B.

To prove C, let $\nu^1$ be a positive function satisfying
$\nu^1 \in C^{\infty}([0,\infty[),$ $\operatorname{supp} \nu^1 \subset
[3/4,\infty[$, (\ref{P01}) and
\begin{equation} \label{h1}
\sum_i (\sesq{\nu^1}{\zcpxd{0}\vol{i}{}})^2 < \infty,
\end{equation}
which is possible since the $\vol{i}{}$ have compact support and by
possibly choosing the $|k_{n+2}|>0$ even smaller. For this given
$\nu^1,$ $\nu\in\mathcal{U}$ is defined as in (P$_1$).

Let $F \in C^{\infty}(\mathbb{R})$ be a positive function satisfying
$\operatorname{supp} F \subset [0,2]$ and $F(1)=1.$
$Y$ is the unique $(\mathcal{F}_t)$-adapted process satisfying
\begin{equation} \label{proof eq 1}
Y_t=1+ \int_0^t F(Y_s)\,dM_s, \qquad t \in\mathbb{T},
\end{equation}
where $M$ is the square integrable $Q$-martingale defined by
\begin{equation} \label{proof eq 2}
M_t=\sum_i \sesq{\nu^1}{\zcpxd{0}\vol{i}{}} \wienerq{i}{t}.
\end{equation}
$X=Y_{\timeh}$ is a positive bounded smooth $\mathcal{F}$ measurable
random variable.

$\mu$ is defined by formulas (\ref{P1.1}) and (\ref{P1}),
with\vspace*{1pt}
$\alpha_t=F(Y_t),$ and it is a generalized integrand. This easily
follows by introducing the sequence $\nu^{(n)} \in H'$ defined by
$\sesq{\nu^{(n)}}{f},$ $f \in H,$ where
\begin{equation} \label{proof eq nu n}
\bigl\langle\nu^{(n)},f\bigr\rangle = \bigl\langle\nu^{(n)1},f\bigr\rangle,\qquad \mbox{if } f \in K
\mbox{ and } \bigl\langle\nu^{(n)},\zcpxd{0}\bigr\rangle=0
\end{equation}
and where $\nu^{(n)1}= \nu^1 g_n$ for a sequence of positive
$C^\infty$ cut-off functions $g_n.$ We here choose $g_n(x)=1$ for $0
\leq x \leq n$ and $g_n(x)=0$ for $x \geq n+1.$ The sequence $\nu^{(n)}
\in H'$ then satisfies (C$_1$) and (C$_2$), which follows similarly as in
the proof of~B.

The decomposition (\ref{eq1 main th1}) of $X$ is valid with $x=1,$ so
$(1, \mu)$ is a generalized hedging portfolio of $X.$

The discounted risk-free investment at $t,$ given by the generalized
portfolio $(1, \mu^{(n)})$ is
\begin{equation} \label{proof eq 3}
a^n_t \zcpxd{t}(0)=Y_t - \bigl\langle\mu^{(n)}_t,\zcpxd{t}\bigr\rangle=Y_t - \alpha
_t \bigl\langle\nu^{(n)1},\zcpxd{0}\bigr\rangle.
\end{equation}
We now choose $\nu^{1}$ and possibly further restrict the $k_i,$ which
is possible, such that
\begin{equation} \label{proof eq 4}
\lim_{n \rightarrow\infty} \bigl\langle\nu^{(n)1},\zcpxd{0}\bigr\rangle=\infty
\end{equation}
and such that the condition in (\ref{h1}) is satisfied. This proves the
part $C=-\infty$ of (P$_2$). The statements for $C$ finite and
$C=+\infty$ are proved so similarly to those in B, that we omit the
proof.
\end{pf*}

\section*{Acknowledgments}
The author thanks Bruno Bouchard for interesting remarks and constructive criticism.


\printaddresses


\begin{thebibliography}{11}

\bibitem{BjKaRu97}
\begin{barticle}[mr]
\bauthor{\bsnm{Bj{\"o}rk},~\bfnm{Tomas}\binits{T.}},
  \bauthor{\bsnm{Kabanov},~\bfnm{Yuri}\binits{Y.}} \AND
  \bauthor{\bsnm{Runggaldier},~\bfnm{Wolfgang}\binits{W.}}
(\byear{1997}).
\btitle{Bond market structure in the presence of marked point processes}.
\bjournal{Math. Finance}
\bvolume{7}
\bpages{211--239}.
\bid{doi={10.1111/1467-9965.00031}, mr={1446647}}
\end{barticle}
\endbibitem

\bibitem{BjMaKaRu97}
\begin{barticle}[vtex]
\bauthor{\bsnm{Bj\"ork},~\bfnm{T.}\binits{T.}},
  \bauthor{\bsnm{Masi},~\bfnm{G.}\binits{G.}},
  \bauthor{\bsnm{Kabanov},~\bfnm{Y.}\binits{Y.}} \AND
  \bauthor{\bsnm{Runggaldier},~\bfnm{W.}\binits{W.}}
  (\byear{1997}).
\btitle{Toward a general theory of bond markets.}
\bjournal{Finance
  Stoch.}
  \bvolume{1}
  \bpages{141--174}.
\end{barticle}
\endbibitem

\bibitem{CarmonaTehr}
\begin{barticle}[mr]
\bauthor{\bsnm{Carmona},~\bfnm{Rene}\binits{R.}} \AND
  \bauthor{\bsnm{Tehranchi},~\bfnm{Michael}\binits{M.}}
(\byear{2004}).
\btitle{A characterization of hedging portfolios for interest rate contingent
  claims}.
\bjournal{Ann. Appl. Probab.}
\bvolume{14}
\bpages{1267--1294}.
\bid{doi={10.1214/105051604000000297}, mr={2071423}}
\end{barticle}
\endbibitem

\bibitem{DaPratoZabczyk}
\begin{bbook}[mr]
\bauthor{\bsnm{Da~Prato},~\bfnm{Giuseppe}\binits{G.}} \AND
  \bauthor{\bsnm{Zabczyk},~\bfnm{Jerzy}\binits{J.}}
(\byear{1992}).
\btitle{Stochastic Equations in Infinite Dimensions}.
\bseries{Encyclopedia of Mathematics and Its Applications}
\bvolume{44}.
\bpublisher{Cambridge Univ. Press}, \baddress{Cambridge}.
\bid{doi={10.1017/CBO9780511666223}, mr={1207136}}
\end{bbook}
\endbibitem

\bibitem{DonnoPratelli2004}
\begin{barticle}[mr]
\bauthor{\bsnm{De~Donno},~\bfnm{Marzia}\binits{M.}} \AND
  \bauthor{\bsnm{Pratelli},~\bfnm{Maurizio}\binits{M.}}
(\byear{2004}).
\btitle{On the use of measure-valued strategies in bond markets}.
\bjournal{Finance Stoch.}
\bvolume{8}
\bpages{87--109}.
\bid{doi={10.1007/s00780-003-0102-7}, mr={2022980}}
\end{barticle}
\endbibitem

\bibitem{DonnoPratelli2005}
\begin{barticle}[mr]
\bauthor{\bsnm{De~Donno},~\bfnm{M.}\binits{M.}} \AND
  \bauthor{\bsnm{Pratelli},~\bfnm{M.}\binits{M.}}
(\byear{2005}).
\btitle{A theory of stochastic integration for bond markets}.
\bjournal{Ann. Appl. Probab.}
\bvolume{15}
\bpages{2773--2791}.
\bid{doi={10.1214/105051605000000548}, mr={2187311}}
\end{barticle}
\endbibitem

\bibitem{IEET2005}
\begin{barticle}[mr]
\bauthor{\bsnm{Ekeland},~\bfnm{Ivar}\binits{I.}} \AND
  \bauthor{\bsnm{Taflin},~\bfnm{Erik}\binits{E.}}
(\byear{2005}).
\btitle{A theory of bond portfolios}.
\bjournal{Ann. Appl. Probab.}
\bvolume{15}
\bpages{1260--1305}.
\bid{doi={10.1214/105051605000000160}, mr={2134104}}
\end{barticle}
\endbibitem

\bibitem{IEET2007}
\begin{bincollection}[vtex]
  \bauthor{\bsnm{Ekeland},~\bfnm{Ivar}\binits{I.}}
  \AND
  \bauthor{\bsnm{Taflin},~\bfnm{Erik}\binits{E.}}
(\byear{2007}).
\btitle{Optimal bond portfolios}. In
\bbooktitle{Paris--{P}rinceton {L}ectures on {M}athematical {F}inance 2004}.
\bseries{Lecture Notes in Math.}
\bvolume{1919}
\bpages{51--102}.
\bpublisher{Springer}, \baddress{Berlin}.
\bid{doi={10.1007/978-3-540-73327-0}, mr={2384671}}
\end{bincollection}
\endbibitem

\bibitem{Pham2003}
\begin{barticle}[vtex]
\bauthor{\bsnm{Pham},~\bfnm{Huy{\^e}n}\binits{H.}}
(\byear{2003}).
\btitle{A predictable decomposition in an infinite assets model with jumps.
  {A}pplication to hedging and optimal investment}.
\bjournal{Stoch. Stoch. Rep.}
\bvolume{75}
\bpages{343--368}.
\bid{doi={10.1080/104511203100001621237}, mr={2017783}}%
\end{barticle}
\endbibitem

\bibitem{Protter2004}
\begin{bbook}[vtex]
\bauthor{\bsnm{Protter},~\bfnm{Philip~E.}\binits{P.~E.}}
(\byear{2004}).
\btitle{Stochastic Integration and Differential Equations},
\bedition{2nd} ed.
\bseries{Applications of Mathematics}
\bvolume{21}.
\bpublisher{Springer}, \baddress{Berlin}.
\bid{mr={2020294}}
\end{bbook}
\endbibitem

\bibitem{ETBondCompleteness}
\begin{barticle}[mr]
\bauthor{\bsnm{Taflin},~\bfnm{Erik}\binits{E.}}
(\byear{2005}).
\btitle{Bond market completeness and attainable contingent claims}.
\bjournal{Finance Stoch.}
\bvolume{9}
\bpages{429--452}.
\bid{doi={10.1007/s00780-005-0156-9}, mr={2211717}}
\end{barticle}
\endbibitem

\end{thebibliography}
\end{document}